\documentclass[oneside,reqno,12pt]{amsart}
\usepackage[greek,english]{babel}
\usepackage{amsthm}
\usepackage{amsbsy}
\usepackage{amsfonts}
\usepackage{graphicx}
\usepackage{hyperref}
\hypersetup{colorlinks=true, citecolor=blue, linkcolor=red}

\newtheorem{thm}{Theorem}
\newtheorem{cor}{Corollary}

\newtheorem{pro}{Proposition}
\newtheorem{rem}{Remark}

\title[Uniqueness for the scalar Dafermos regularization]{Uniqueness for the Dafermos regularization viscous wave fan profiles for
Riemann solutions of scalar hyperbolic conservation laws}
\author{Christos Sourdis}
\address{General Lyceum of Moires, Heraklion, Crete, Greece.}
\email{sourdis@uoc.gr}

\begin{document}

\maketitle
\begin{abstract}
We prove the uniqueness of solutions to the Dafermos regularization viscous wave fan profiles for
Riemann solutions of scalar hyperbolic conservation laws. We emphasize that our results are not restricted to the small self-similar viscosity regime. We rely on suitable adaptations of  Serrin's sweeping principle and the sliding method from the qualitative theory of semilinear elliptic PDEs. In order to illustrate the delicacy of our result, we prove the existence of an unbounded solution in the case of Burgers equation. Lastly, we can combine aspects of  these results in order to give a precise description of the Dafermos regularization of rarefaction waves of Burgers equation.
\end{abstract}


\section{Introduction}
The  Riemann problem for a single conservation law in a single space variable is to determine a self-similar (generally weak) solution $U$ of
 \begin{equation}\label{eqEq11}
U_t+f(U)_x=0,\ \ x \in \mathbb{R}, \ t>0,
\end{equation}
  with the initial condition
\begin{equation}\label{eqBC}
U(x,0)=\left\{\begin{array}{cc}
                u_L & \textrm{for}\ x<0, \\
                u_R & \textrm{for}\ x>0,
              \end{array}
 \right.
\end{equation}
where $f\in C^1$ (see for instance \cite[Ch. IX]{dafbook} or \cite[Ch. 16]{smoller}).
A solution is self-similar if
\begin{equation}\label{eqscalinas}U(x,t)=u(\xi),\  \xi=x/t.\end{equation}

 The corresponding Dafermos regularization \cite{dafclassic, ka, tupciev} is\begin{equation}\label{eqdaf}
                                                                          U_t+f(U)_x=\varepsilon t U_{xx},\ \ x \in \mathbb{R}, \ t>0,\end{equation}
together with (\ref{eqBC}),
 where $\varepsilon>0$ (typically small).
Self-similar solutions of (\ref{eqdaf}) with respect to the scaling in (\ref{eqscalinas}) represent viscous
wave fan profiles for Riemann solutions.   The corresponding ODE  problem for $u$ is
\begin{equation}\label{eqEq}
\varepsilon u_{\xi\xi}=\left(f'(u)-\xi\right)u_\xi,\ \xi\in\mathbb{R};\ u(-\infty)=u_L,\ u(+\infty)=u_R.
\end{equation} Throughout this paper, solutions to the above ODE are understood in the classical sense (i.e. $u\in C^2(\mathbb{R})$).

Existence of a solution to (\ref{eqEq}) is known to hold  for any $\varepsilon>0$ by a fixed point argument (see \cite{dafclassic,tzavaras}) and such solutions are monotone. In fact, their  monotonicity  (which follows at once from the uniqueness of initial value problems for ODEs) provides the uniform a-priori estimates that are needed in the  existence proof. We also refer to \cite{schecter} for the construction of solutions for small $\varepsilon>0$ using geometric singular perturbation theory.
On the other hand, to the best of our knowledge, a uniqueness result for  (\ref{eqEq}) seems to be missing (see the many related references for Section 9.8  of \cite{dafbook}). In fact, uniqueness seems to be unknown even for the Dafermos regularization of Burgers equation (where $f(u)=u^2/2$) with $\varepsilon>0$ small, see  Remark 3.2 in \cite{liu}.   In this regard, we will establish in this note  the following uniqueness result for (\ref{eqEq}).

\begin{thm}\label{thmMain}
If   $f'$ is  Lipscitz continuous on $[\min\{u_L,u_R\},\max\{u_L,u_R\}]$ if $u_L\neq u_R$ or $f\in C^1(\mathbb{R})$ if $u_L=u_R$, and $\varepsilon>0$, then there exists a unique solution to (\ref{eqEq}).
\end{thm}

 A Riemann solution $U(x, t) = u^*(\xi)$ of (\ref{eqEq11}) and (\ref{eqBC}) satisfies the \emph{viscous wave
fan profile criterion} if and only if there exists a solution $u_\varepsilon(\xi)$ of (\ref{eqEq})  such that
\[u_\varepsilon(\xi) \to  u^*(\xi)\ \textrm{as}\ \varepsilon\to 0\ \textrm{in}\ L^1_{loc}(\mathbb{R}).\]
As a direct consequence of Theorem \ref{thmMain}, we get the following.
\begin{cor}\label{cor}
If   $f'$ is  Lipscitz continuous on $[u_L,u_R]$, then there exists at most one  solution to (\ref{eqEq11}) and (\ref{eqBC}) which satisfies the viscous wave
fan profile criterion.
\end{cor}

In particular, Theorem \ref{thmMain} shows the uniqueness of the  composite-wave solutions that were constructed in \cite{schecter} for small $\varepsilon>0$ via geometric singular perturbation theory.
Moreover, we can provide an answer to
 Remark 3.2 in \cite{liu} concerning the uniqueness of solutions to the corresponding Dafermos regularization for Burgers equation, i.e., the ODE problem
 \begin{equation}\label{eqEqBurger}
\varepsilon u_{\xi\xi}=(u-\xi)u_\xi,\ u(-\infty)=u_L,\ u(+\infty)=u_R.
\end{equation}
In the aforementioned paper \cite{liu},
it was proven that  the three-dimensional dynamical
system associated to (\ref{eqEqBurger}) possesses    a  first integral (see also Remark \ref{remliu} herein), and an explicit formula for a normally
hyperbolic invariant curve  was given, as well as
for its stable and unstable manifolds. Then, using these results,   a direct, global,
and rather detailed treatment of the singularly perturbed problem (\ref{eqEqBurger}) was given for solutions which converge as $\varepsilon\to 0$ to rarefaction  and shock waves for $u_L<u_R$ and $u_R<u_L$ respectively (if $u_L=u_R$, then clearly there is only the constant solution).
We refer the interested reader to \cite{slemrod} for a comparison of solutions to (\ref{eqEqBurger}) with solutions of the standard viscous regularization
 \[
U_t+\left(\frac{U^2}{2}\right)_x=\varepsilon U_{xx},\ \ x \in \mathbb{R}, \ t>0,
\]
as $\varepsilon\to 0$.

In order to highlight the delicacy of Theorem \ref{thmMain}, we note that uniqueness may not hold without both 'boundary conditions' in (\ref{eqEq}). Indeed, we have  the following proposition which is of independent interest and will be of use in the subsequent Proposition \ref{procorner} for the study of the singular perturbation problem (\ref{eqEqBurger}).
\begin{pro}\label{proburgers}
There exists a unique, increasing solution to the ODE in (\ref{eqEqBurger}) with $\varepsilon=1$, i.e.
 \begin{equation}\label{eqEqBurger1111}
  u_{\xi\xi}=(u-\xi)u_\xi, \ \xi\in \mathbb{R},
\end{equation}
such that
\begin{equation}\label{eqAsy}
   \max\{0,\xi\}<u(\xi)<\max\{0,\xi\}+Ce^{-c|\xi|}\ \textrm{in}\ \mathbb{R},
\end{equation}
for some constants $c,C>0$. Actually, the constant $c$ in the exponent can be taken equal to $1$ (see Remark \ref{remliu} for more details).
\end{pro}

Armed with the above proposition, we can employ once more the sliding method in order to establish the following result which refines \cite[Thm. 3.1]{liu}.
\begin{pro}\label{procorner}
  The solution of (\ref{eqEqBurger}) with $u_L<u_R$ satisfies
  \[
  u_\varepsilon(\xi)=\sqrt{\varepsilon}\texttt{U}\left(\frac{\xi-u_L}{\sqrt{\varepsilon}}\right)+u_L+\sqrt{\varepsilon} O\left(e^{-1/\sqrt{\varepsilon}}\right),\ \xi\leq \frac{u_L+u_R}{2},
  \]
  and
  \begin{equation}\label{eqsyms}
    u_\varepsilon(u_L+u_R-\xi)+u_\varepsilon(\xi)=u_L+u_R,\ \xi\in \mathbb{R},
  \end{equation}
  where $\texttt{U}$ is as in Proposition \ref{proburgers}, and the (normalized) remainder $e^{1/\sqrt{\varepsilon}}O\left(e^{-1/\sqrt{\varepsilon}}\right)$ is uniformly bounded in $\mathbb{R}$ as $\varepsilon\to 0$.
\end{pro}

The solutions provided by the above proposition converge to \emph{rarefaction waves} for (\ref{eqEq11})-(\ref{eqBC}) as $\varepsilon\to 0$. On the other hand, if $u_L>u_R$ then (\ref{eqEq11})-(\ref{eqBC}) admits the following \emph{shock wave} solution:
\[u^*(\xi)=\left\{\begin{array}{ll}
             u_L, & \xi\leq (u_L+u_R)/2, \\
             u_R, & \xi> (u_L+u_R)/2.
           \end{array}\right.
\]The corresponding solutions to (\ref{eqEq}) can be constructed using well known arguments of geometric singular perturbation theory (see for example \cite{ku}). We point out that in that case there is no loss of normal hyperbolicity as opposed to the case in Proposition \ref{procorner} (see \cite{liu}).

Our method of proof of Theorem \ref{thmMain}  is based on Serrin's sweeping principle (see \cite{pucci,serg}) and the famous sliding method (see for instance \cite{bereSliding,bhm} and the references therein) from the qualitative theory of elliptic PDEs. In the same spirit, the proof of Proposition \ref{proburgers} is based on the method of super and subsolutions.
Finally, the proof of Proposition \ref{procorner} is highly motivated by the sliding argument in the proof of Theorem \ref{thmMain}.

The rest of the paper is devoted to the proofs of the above results and to some related remarks.
\section{Proofs of the main results}\label{sec2}
\subsection*{Proof of Theorem \ref{thmMain}}
\begin{proof}
As we have already mentioned, existence of a solution to (\ref{eqEq}) is known to hold for any $\varepsilon>0$. Moreover, by the classical uniqueness result for the initial value problem associated to the first order linear ODE for $u_\xi$, it follows that solutions of (\ref{eqEq}) are strictly monotone unless they are identically constant. 
We will distinguish the following three cases.

\underline{\textbf{$u_L=u_R.$}} As we have just mentioned, in that case the only solution is the constant one.

\underline{\textbf{$u_L<u_R.$}} Let $u,v$ be two solutions of (\ref{eqEq}). We will adapt the famous sliding method  (see for instance \cite[Thm. 1]{bhm}) in order to show that $u\geq v$.  Once this is established,  the desired uniqueness property will follow simply by exchanging the roles of $u$ and $v$. We note that in the problem at hand, in contrast to \cite{bhm} and related references where the sliding method is employed, the aforementioned  monotonicity property
\begin{equation}\label{eqMono}
u_\xi,\ v_\xi>0,\ \xi\in \mathbb{R},
\end{equation}
will play an important role throughout the proof.

Let us consider the translations
\begin{equation}\label{equl2}
   u_\lambda(\xi)=u(\xi+\lambda), \ \ \lambda\geq 0.
\end{equation}
We note that, thanks to (\ref{eqMono}), $u_\lambda$ is strictly increasing with respect to $\lambda\geq 0$.
The main observation is that, for $\lambda>0$, we have
\begin{equation}\label{eqEq2}
\varepsilon ( u_\lambda)_{\xi\xi}<\left( f'(u_\lambda)-\xi\right)( u_\lambda)_\xi,\ \xi\in \mathbb{R};\  u_\lambda(-\infty)=u_L,\  u_\lambda(+\infty)=u_R,
\end{equation}
 i.e., $u_\lambda$ is a strict supersolution of (\ref{eqEq}).
 Indeed, setting $\xi+\lambda$ with $\lambda>0$ in place of $\xi$ in (\ref{eqEq}) gives
\[
\varepsilon u_{\xi\xi}(\xi+\lambda)=\left(f'\left(u(\xi+\lambda)\right)-\xi-\lambda\right)u_\xi(\xi+\lambda)\stackrel{(\ref{eqMono})}{<}\left(f'\left(u(\xi+\lambda)\right)-\xi\right)u_\xi(\xi+\lambda),
\]
and   the desired relation  (\ref{eqEq2}) follows readily.

To start  the sliding process, we will show that  there exists a large $\Lambda>0$   such that
\begin{equation}\label{eqorder0}
u_\Lambda>v\ \textrm{in}\ \mathbb{R}.
\end{equation}
To this end, we first note that, given any $M>0$, by the asymptotic behaviour of $u$ from (\ref{eqEq}) there clearly exists a $\Lambda(M)>0$ sufficiently large such that
\begin{equation}\label{eqM}
  u_\lambda>v,\ \xi\in [-M,M],\ \textrm{if}\ \lambda\geq \Lambda(M).
\end{equation}
Then, using (\ref{eqEq}) for $v$   and (\ref{eqEq2}), we will infer that the above strict inequality continues to hold for $|\xi|>M$, provided that $M$ is chosen sufficiently large.
Indeed, from the aforementioned relations, we see that
\[\psi=u_\lambda-v\]
satisfies
\[
\varepsilon \psi_{\xi\xi}<\left(f'(u_{\lambda})-\xi\right)\psi_\xi+v_\xi\left(f'(u_{\lambda})-f'(v)\right),\ \xi\in \mathbb{R}.
\]
We can equivalently rewrite the above relation as
\begin{equation}\label{eqpsi22222}L(\psi)=
\varepsilon \psi_{\xi\xi}-\left(f'(u_{\lambda})-\xi\right)\psi_\xi-v_\xi Q(\xi)\psi<0,\ \xi\in \mathbb{R},
\end{equation}
where
\[
Q(\xi)=\left\{\begin{array}{ll}
                \frac{f'\left(u_{\lambda}(\xi)\right)-f'\left(v(\xi)\right)}{u_\lambda(\xi)-v(\xi)} & \textrm{if}\ u_\lambda(\xi)\neq v(\xi),  \\
                0 & \textrm{if}\ u_\lambda(\xi)= v(\xi).
              \end{array}
 \right.
\]
Since $f'$ is  Lipschitz continuous on $[u_L,u_R]$, we have that $Q$ is bounded in absolute value by the Lipschitz constant of $f'$ over $[u_L,u_R]$.
Summarizing, we have
\begin{equation}\label{eqsumar}
  L(\psi)<0\ \textrm{if}\ |\xi|>M;\ \psi(\pm M)>0,\ \psi\to 0 \ \textrm{super-exponentially fast as}\ |\xi|\to \infty,
\end{equation}
(for the last property we refer to (3.9) in \cite{dafclassic}). We will show that if $M>0$ is sufficiently large then \begin{equation}\label{eqout}\psi>0\ \textrm{for}\ |\xi|>M\end{equation} by a maximum principle type argument, despite of the fact that $Q$ may be sign-changing. To this end, by (1.6) in \cite{bnv}, it suffices to prove that there exists a   $g\in C^2$ such that
\begin{equation}\label{eqLg}
  L(g)\leq 0\ \textrm{if}\ |\xi|>M,\ g>0\ \textrm{if}\ |\xi|\geq M \ \textrm{and}\ \lim_{|\xi|\to \infty}\frac{\psi}{g}=0.
\end{equation}
It turns out that \[g(\xi)=e^{-|\xi|}\] satisfies (\ref{eqLg}), provided that $M>0$ is chosen sufficiently large. Indeed, by virtue of the asymptotic behaviour of $\psi$ from
(\ref{eqsumar}), it remains to verify the first property in (\ref{eqLg}). We will do this only for $\xi>M$ as the case $\xi<-M$ can be treated similarly. A simple calculation shows that
\[
L(g)=\left(\varepsilon + f'(u_{\lambda})-\xi-v_\xi Q(\xi) \right)e^{-\xi}<0\ \textrm{if}\ \xi>M,
\]
provided that $M$ is chosen sufficiently large (since $f'(u_\lambda), v_\xi, Q$ are bounded). In fact, let us from now on choose
\[M=
1+\varepsilon + \|f'\|_{L^\infty(u_L,u_R)}+ \|v_\xi \|_{L^\infty(\mathbb{R})}\|f'\|_{C^{0,1}(u_L,u_R)},
\]
which clearly satisfies the above properties.
The desired relation (\ref{eqorder0}) follows at once by combining (\ref{eqM}) and (\ref{eqout}).

We can now define $\lambda_0\geq 0$ by
\begin{equation}\label{eqlam0}
  \lambda_0=\inf\left\{\lambda \in (0,\Lambda]\ :\ u_\lambda\geq v \ \textrm{in}\ \mathbb{R} \right\}.
\end{equation}
Obviously, by continuity, we have
\begin{equation}\label{eqbarber}
u_{\lambda_0} \geq v \ \textrm{in}\ \mathbb{R}.
\end{equation}
 We will establish that $\lambda_0=0$. To this end, arguing by contradiction, let us suppose that $\lambda_0>0$.
 If \[\min_{[-M,M]}(u_{\lambda_0}-v)=0,\]
  in light of (\ref{eqbarber}), we can easily reach a contradiction by applying the strong maximum principle (see for instance \cite[Thm. 2.8.4]{pucci}) in the strict differential inequality (\ref{eqpsi22222}) (with $\lambda=\lambda_0$).
 It remains to consider the case where\[\min_{[-M,M]}(u_{\lambda_0}-v)>0.\]
 However, since $u_\lambda$ is continuous with respect to $\lambda$,  the above relation implies that there exists some small $\epsilon\in (0,\lambda_0)$ such that
\[\min_{[-M,M]}(u_{\lambda_0-\epsilon}-v)>0.\]
In turn, by the maximum principle as before    we get
\[
u_{\lambda_0-\epsilon}>v\ \textrm{in}\ \mathbb{R},
\]
which contradicts the definition of $\lambda_0$. We thus conclude that $\lambda_0=0$. Hence, we infer from (\ref{eqbarber}) that $v\leq u$ as desired.

\underline{\textbf{$u_L>u_R.$}}
The proof is analogous to the previous case, in the sense that it is based on a continuity argument and the maximum principle. However, here solutions satisfy \begin{equation}\label{eqMono1}
u_\xi<0,\ \xi\in \mathbb{R},
\end{equation}
 and we cannot use the previous sliding argument since we have the 'wrong' inequality in (\ref{eqEq2}). Instead, we will employ Serrin's sweeping principle
  (see \cite[Thm. 2.3.5]{pucci} or \cite{serg}). More precisely, we will 'sweep' with the family
\begin{equation}\label{equl}
u_\lambda(\xi)=u(\xi-2K\lambda)+\lambda,\ \ \lambda\geq0,
\end{equation}
where $K>0$ denotes the Lipschitz constant of $f'$, i.e.
\begin{equation}\label{eqlipschitz}
|f'(u_1)-f'(u_2)|\leq K|u_1-u_2|,\ u_1,u_2\in [u_R,u_L].
\end{equation} We note in passing that (\ref{equl}) is motivated by Remark \ref{rem2222} below for Burgers equation.

We claim that $u_\lambda$ is a strict supersolution to (\ref{eqEq}) for $\lambda>0$, that is
  \begin{equation}\label{eqEq3}
\varepsilon ( u_\lambda)_{\xi\xi}<\left( f'(u_\lambda)-\xi\right)( u_\lambda)_\xi,\ \xi\in \mathbb{R};\  u_\lambda(-\infty)=u_L+\lambda,\  u_\lambda(+\infty)=u_R+\lambda.
\end{equation} Indeed, setting $\xi-2K\lambda$ with $\lambda>0$ in place of $\xi$ in (\ref{eqEq}), we find
\[\begin{array}{ccl}
  \varepsilon  u_{\xi\xi}(\xi-2K\lambda) & = & \left[f'\left(u(\xi-2K\lambda) \right)+2K\lambda-\xi \right]u_{\xi}(\xi-2K\lambda) \\
    \ \ \\
    \textrm{via} \ (\ref{eqMono1}), (\ref{eqlipschitz}):  & \leq & \left[f'\left(u(\xi-2K\lambda)+\lambda \right)-K\lambda+2K\lambda-\xi \right]u_{\xi}(\xi-2K\lambda) \\
    \ \ \\
    \textrm{again thanks to} \ (\ref{eqMono1}):& < & \left[f'\left(u(\xi-2K\lambda)+\lambda \right)-\xi \right]u_{\xi}(\xi-2K\lambda), \\
  \end{array}
 \]and the desired relation (\ref{eqEq3}) follows at once.

 Armed with the above information, the proof proceeds along the same lines as in the previous case. We just point out that, thanks to  the  asymptotic behaviour of $u_\lambda$ from (\ref{eqEq3}), the continuity argument is actually simpler and there is no need for a splitting of the real line here. We omit the details.
\end{proof}

\begin{rem}\label{remconvex}
  It is easy to see that the proof of Theorem \ref{thmMain} in the case where  $u_L<u_R$ applies under the weaker assumption that $f \in C^1(u_L,u_R)$, provided that $f$ is convex
  near $u_L$ and $u_R$. Indeed, the last assumption implies that $Q(\xi)\geq 0$ if $|\xi|$ is sufficiently large, and one can plainly take $g\equiv 1$.
  
  If   $u_L> u_R$ and $f\in C^1(u_R,u_L)$   is concave in $(u_R,u_L)$, then one can show that there is uniqueness as in the last part of the proof of Theorem \ref{thmMain} by
  noting that
    \[
  u_\lambda=u+\lambda, \ \lambda\geq 0,
  \]
   is a family of supersolutions to (\ref{eqEq}).
     \end{rem}

\subsection*{Proof of Proposition \ref{proburgers}}
\begin{proof}
  The function \[\underline{u}(\xi)=\max\{0,\xi\}\] provides a weak subsolution to (\ref{eqEqBurger1111}) (see for example \cite{berestlions}).
Indeed, $\underline{u}$ satisfies (\ref{eqEqBurger1111}) for $\xi\neq 0$ and $\underline{u}_\xi(0^-)\leq \underline{u}_\xi(0^+)$.

We claim that
\[
\bar{u}(\xi)=\left\{\begin{array}{ll}
                      \int_{-\infty}^{\xi}e^{-\frac{t^2}{2}}dt, & \xi\leq 0, \\
                        &   \\
                      \xi+\int_{-\infty}^{0}e^{-\frac{t^2}{2}}dt, & 0< \xi\leq 1, \\
                        &   \\
                      \xi+\left(\int_{-\infty}^{0}e^{-\frac{t^2}{2}}dt\right) e^{-\frac{\xi-1}{L}},& \xi>1,
                    \end{array}
 \right.
\]
is a weak supersolution to (\ref{eqEqBurger1111}), provided that $L>0$ is chosen sufficiently large. Indeed, letting $I=\int_{-\infty}^{0}e^{- {t^2}/{2}}dt$,  we have
\[
-\bar{u}_{\xi\xi}+(\bar{u}-\xi)\bar{u}_{\xi}=\xi e^{-\frac{\xi^2}{2}}+(\bar{u}-\xi)e^{-\frac{\xi^2}{2}}=\bar{u}e^{-\frac{\xi^2}{2}}>0,\ \xi<0;
\]
\[
-\bar{u}_{\xi\xi}+(\bar{u}-\xi)\bar{u}_{\xi}=\int_{-\infty}^{0}e^{-\frac{t^2}{2}}dt>0, \ 0<\xi<1;
\]
\[\begin{array}{rcl}
    -\bar{u}_{\xi\xi}+(\bar{u}-\xi)\bar{u}_{\xi} & = & -\frac{I}{L^2}e^{-\frac{\xi-1}{L}}
+I e^{-\frac{\xi-1}{L}}\left(1-\frac{I}{L} e^{-\frac{\xi-1}{L}} \right) \\
     &  &  \\
     & = & Ie^{-\frac{\xi-1}{L}}\left(-\frac{1}{L^2}+1-\frac{I}{L} e^{-\frac{\xi-1}{L}} \right)>0,
  \end{array}
\]
$\xi>1$, provided that $L$ is chosen sufficiently large. Moreover, note that
\[
\bar{u}_{\xi}(0^-)=\bar{u}_{\xi}(0^+)\ \textrm{and}\ \bar{u}_{\xi}(1^-)\geq \bar{u}_{\xi}(1^+).
\]
The  claim  that $\bar{u}$ is a weak supersolution to (\ref{eqEqBurger1111}) now follows directly from \cite{berestlions}.

Since $\underline{u}<\bar{u}$, the existence of a solution to (\ref{eqEqBurger1111}) such that $\underline{u}<u<\bar{u}$ follows from \cite{ako,berestlions,serg} (these references consider the case of a bounded domain, nevertheless we can conclude by a standard limiting argument). Clearly, this solution satisfies all the desired properties. We point out that its uniqueness can be shown by a sliding argument as in Theorem \ref{thmMain}.

\end{proof}
\begin{rem}\label{remliu}
It was observed in \cite{liu} that (\ref{eqEqBurger1111}) has the following first integral:
\[
H=\frac{1}{2}(u-\xi)^2-(u_\xi-1)+\ln |u_\xi|.
\]
Therefore, the solution $\texttt{U}$ provided by Proposition \ref{proburgers} satisfies
\begin{equation}\label{eqH0}
  \frac{1}{2}(u-\xi)^2=(u_\xi-1)-\ln u_\xi,\ \xi \in \mathbb{R}.
\end{equation}
By Taylor's expansion, we have
\[
v-1-\ln v=\frac{(v-1)^2}{2} +O\left((v-1)^3\right) \ \textrm{as}\ v \to 1,
\]
where we have employed once more Landau's notation for the remainder.
Hence, we deduce from (\ref{eqH0}) that
\[
\lim_{\xi\to +\infty}\frac{u-\xi}{u_\xi-1}=-1,
\](keep in mind that $0<U_\xi<1$ by the convexity of $U$).
In turn, by L'Hopital's rule we obtain that
\[
\lim_{\xi\to +\infty}\frac{\ln(u-\xi)}{\xi}=-1.
\]
We thus conclude that the constant $c$ in the exponent in (\ref{eqAsy}) can  actually be taken equal to $1$ as $\xi\to +\infty$. On the other side, it follows from the proof of Proposition \ref{proburgers} that the corresponding decay rate as $\xi\to -\infty$ is super-exponentially fast.
\end{rem}
\begin{rem}\label{rem2222}
Assume that $u$ satisfies the ODE in (\ref{eqEqBurger}). Then   \[
u_\lambda(\xi)=u(\xi-\lambda)+\lambda,\ \ \lambda \in \mathbb{R},
\]
   solves the same ODE. In other words, the ODE in (\ref{eqEqBurger}) is invariant under the above transformation. To see this, note that
setting $\xi-\lambda$ in place of $\xi$ in (\ref{eqEqBurger}) yields
\[
\varepsilon u_{\xi\xi}(\xi-\lambda)=\left(u(\xi-\lambda)+\lambda-\xi\right)u_\xi(\xi-\lambda),
\]
and the desired conclusion follows at once.
\end{rem}
\subsection*{Proof of Proposition \ref{procorner}}
\begin{proof}
  It follows from the uniqueness property in Theorem \ref{thmMain} and \cite[Prop. 2.7]{liu} that the solution to (\ref{eqEqBurger}) satisfies the odd symmetry property (\ref{eqsyms}). Therefore, problem (\ref{eqEqBurger}) reduces to
   \begin{equation}\label{eqEqBurgerSym}
\varepsilon u_{\xi\xi}=(u-\xi)u_\xi,\ \xi<\frac{u_L+u_R}{2};\ u(-\infty)=u_L,\ u\left(\frac{u_L+u_R}{2}\right)=\frac{u_L+u_R}{2}.
\end{equation}

It is easy to see that
\begin{equation}\label{eq1}
v_\varepsilon(\xi)=\sqrt{\varepsilon}\texttt{U}\left(\frac{\xi-u_L}{\sqrt{\varepsilon}}\right)+u_L
\end{equation}
satisfies
\[
  \varepsilon v_{\xi\xi}=(v-\xi)v_\xi,\ \xi<\frac{u_L+u_R}{2}\ \textrm{(see also Remark \ref{rem2222})};\]\[ v(-\infty)=u_L,\ v\left(\frac{u_L+u_R}{2}\right)=\frac{u_L+u_R}{2}+\sqrt{\varepsilon}O\left(e^{-1/\sqrt{\varepsilon}}\right),\ \textrm{as}\ \varepsilon\to 0.
\]

As in the proof of Theorem \ref{thmMain}, since $(v_{\varepsilon,\lambda})_\xi>0$, we find that
\begin{equation}\label{eq2}
  v_{\varepsilon,\lambda}(\xi)=v_\varepsilon(\xi+\lambda)
\end{equation}
is a (strict)  subsolution to the ODE in (\ref{eqEqBurgerSym}) if $\lambda<0$ and a (strict) supersolution if $\lambda>0$. Furthermore, $v_{\varepsilon,\lambda}$ is strictly increasing with respect to $\lambda$. Moreover, thanks to (\ref{eqAsy}), we have
\[v_{\varepsilon,\lambda}(-\infty)=u_L,\ v_{\varepsilon,\lambda}\left(\frac{u_L+u_R}{2}\right)=\frac{u_L+u_R}{2}+\lambda+\sqrt{\varepsilon}O\left(e^{-1/\sqrt{\varepsilon}}\right)\ \textrm{as}\ \varepsilon\to 0.\] Hence, as in Theorem \ref{thmMain} or Proposition \ref{procorner}, the solution of (\ref{eqEqBurgerSym}) must satisfy
\[
v_{\varepsilon,-\lambda_\varepsilon}<u<v_{\varepsilon,\lambda_\varepsilon},\ \xi\in \mathbb{R},\  \textrm{where}  \ 0<\lambda_\varepsilon<\sqrt{\varepsilon}O\left(e^{-1/\sqrt{\varepsilon}}\right)\ \textrm{as}\ \varepsilon\to 0.
\]
The assertion of the proposition now follows readily from (\ref{eq1}), (\ref{eq2}) and the fact that $\texttt{U}_\xi$ is bounded (in fact,  $0<\texttt{U}_\xi<1$ holds due to the convexity of $\texttt{U}$).
\end{proof}

\end{document}